\documentclass[reqno,10pt]{amsart}

\usepackage{amsmath}
\usepackage{amssymb}
\usepackage{amsthm}
\usepackage{xspace}
\usepackage{pictexwd}
\usepackage{float}  \restylefloat{figure} %\restylefloat{table}
\usepackage[OT2,OT1]{fontenc}
\usepackage{amscd}
\usepackage{enumerate}
\usepackage{graphicx}

\theoremstyle{plain}
\newtheorem{thm}{Theorem}[section]
\newtheorem*{thm*}{Theorem}

\theoremstyle{definition}
\newtheorem{mydef}{Definition}[section]

\theoremstyle{remark}

\def\forcehmode{\hskip0pt\relax}

\let\myskip=\medskip

\def\definebb#1=#2.{\def#1{{{\mathbb #2}^{\vphantom{x}}}}}
\definebb \c=C.
\definebb \r=R.
\definebb \s=S.
\definebb \z=Z.
\definebb \h=H.
\def\Cal#1{\mathcal{#1}}

\def\calD{{\Cal D}}
\def\calM{{\Cal M}}
\def\calT{{\Cal T}}
\def\calV{{\Cal V}}

\def\hyp{\h}

\def\PSL{\MathOpPSL(2,\r)}

\def\dd{\partial}

\def\al{{\alpha}}
\def\be{{\beta}}
\def\Ga{{\Gamma}}

\def\ve{{\varepsilon}}
\def\la{{\lambda}}
\def\si{{\sigma}}

\def\st{\,\,\big|\,\,}
\def\<{\langle}
\def\>{\rangle}

\let\le=\leqslant

\DeclareMathOperator{\Aut}{Aut}

\DeclareMathOperator{\MathOpPSL}{PSL}
\DeclareMathOperator{\Mod}{Mod}

\let\Im=\undefined \DeclareMathOperator{\Im}{Im}
 
\let\mod=\undefined \DeclareMathOperator{\mod}{mod}

\begin{document}

\author[Sergey Natanzon]{Sergey Natanzon}
\address{Faculty of Mathematics, HSE University, Usacheva Str 6, 119048 Moscow, Russian Federation}
%\address{National Research University Higher School of Economics (HSE), Vavilova Street 7, 117312 Moscow, Russia}
%%\address{Laboratory of Quantum Topology, Chelyabinsk State University, Chelyabinsk, Russia}
%%\address{Belozersky Institute of Phys.-Chem. Biology, Moscow State University, Korp.~A, Leninske Gory, 11899 Moscow, Russia}
%\address{Institute of Theoretical and Experimental Physics (ITEP), Moscow, Russia}
%\address{Independent University of Moscow, Bolshoi Vlasevsky Pereulok 11, 119002 Moscow, Russia}
\email{natanzons@mail.ru}
%\email{natanzon@mccme.ru}
\author[Anna Pratoussevitch]{Anna Pratoussevitch}
\address{Department of Mathematical Sciences\\ University of Liverpool\\ Liverpool L69~7ZL, United Kingdom}
\email{annap@liv.ac.uk}

\title{Hyperbolic Groups and Non-Compact Real Algebraic Curves}

\begin{date}  {\today} \end{date}

\thanks{Grant support for S.N.: This work was supported in part by the grant RFBR-20-01-00579.}

%\thanks{Grant support for S.N.: The article was prepared within the framework of the Academic Fund Program
%at the National Research University Higher School of Economics (HSE) in 2015--16 (grant Nr 15-01-0052)
%and supported within the framework of a subsidy granted to the HSE by the Government of the Russian Federation for the implementation of the Global Competitiveness Program.
%The work was supported in part by the grant RFBR-16-01-00409.
%Grant support for A.P.: The work was supported in part by the Leverhulme Trust grant RPG-057.}

\begin{abstract}
In this paper we study the spaces of non-compact real algebraic curves, i.e.\ pairs $(P,\tau)$,
where $P$ is a compact Riemann surface with a finite number of holes and punctures
and $\tau:P\to P$ is an anti-holomorphic involution.
We describe the uniformisation of non-compact real algebraic curves by Fuchsian groups.
We construct the spaces of non-compact real algebraic curves
and describe their connected components.
We prove that any connected component is homeomorphic to a quotient
of a finite-dimensional real vector space by a discrete group
and determine the dimensions of these vector spaces.
\end{abstract}


% 30F50: Klein surfaces
% 30F35: Fuchsian groups and automorphic functions
% 30F60: Teichmueller theory

%\keywords{Higher spin bundles, real forms, Riemann surfaces, Klein surfaces, Arf functions, lifts of Fuchsian groups}

\maketitle

%\tableofcontents

\bigskip
\hfill\emph{In memory of Ernest Borisovich Vinberg}

\bigskip
\section{Introduction}

\myskip
It is well known that the category of complex algebraic curves is isomorphic to the category of compact Riemann surfaces.
For a complex algebraic curve generated by the polynomial~$F(x,y)$ with complex coefficients, 
the corresponding compact Riemann surface~$P$ is obtained as the regularisation of the set of complex solutions of the equation~$F(x,y)=0$.
If the Riemann surface~$P$ has a finite number of holes and punctures, 
we will say that the corresponding complex curve has holes and punctures.
Such complex curves play an important role in current research, see~\cite{Wi}.

\myskip
Similarly the category of real algebraic curves is isomorphic to the category of pairs~$(P,\tau)$,
where $P$ is a compact Riemann surface and $\tau:P\to P$ is an anti-holomorphic involution on~$P$.
For a real algebraic curve generated by the polynomial~$F(x,y)$ with real coefficients, 
the corresponding pair~$(P,\tau)$ consists of the compact Riemann surface~$P$ as above 
and the involution~$\tau$ which is generated by the complex conjugation $(x,y)\mapsto(\bar x,\bar y)$.
The set of fixed points of the involution~$\tau$ is called the {\it set of real points\/} of the real curve~$(P,\tau)$
and is denoted by~$P^{\tau}$.
If the Riemann surface~$P$ has a finite number of holes and punctures,
we will say that $(P,\tau)$ is a real curve with holes and punctures.
Those holes and punctures that are invariant under the involution~$\tau$
will be called real holes and real punctures.
The remaining holes and punctures occur in pairs which are mapped by the involution~$\tau$.

\myskip
The set of real points of a compact real curve decomposes into pairwise disjoint simple closed smooth curves
that are called {\it ovals\/}.
The {\it topological type\/} of the compact real curve $(P,\tau)$ is determined by the triple $(g,k,\ve)$,
where $g$ is the genus of~$P$,
$k$ is the number of ovals
and $\ve\in\{0,1\}$ with $\ve=1$ if~$P\backslash P^{\tau}$ is orientable and $\ve=0$ otherwise.

\myskip
We can compactify a real curve~$(P,\tau)$ by closing every hole with a disk and every puncture with a point
and extending~$\tau$ to an anti-holomorphic involution~$\hat\tau$ on the resulting surface~$\hat P$.
The ovals of~$(\hat P,\hat\tau)$ are called {\it compactified ovals\/}.
The {\it type\/} of a compactified oval is given by the cyclic sequence of glued in disks and points.
The {\it topological type\/} of the real curve $(P,\tau)$ is determined by the type~$(g,k,\ve)$ of the corresponding compactified real curve $(\hat P,\hat\tau)$,
the numbers~$2n_I$ and~$2m_I$ of non-real holes and punctures
and the types of compactified ovals.
The numbers~$n_R$ and~$m_R$ of real holes and punctures respectively
are determined by the type of the real curve.
In this paper we will assume that $2g+2n_I+n_R+2m_I+m_R>2$
as the other cases are easier, but require different techniques.

\myskip
The aim of this paper is to prove the following description of the moduli space of real curves with holes and punctures.

\begin{thm}
\label{thm-moduli}
The moduli space~$\calM_t$ of all real algebraic curves of type
\[t=(g,k,\ve|2n_I,2m_I,\text{types of compactified ovals},n_R,m_R)\]
is not empty if and only if
%the type~$(g,k,\ve)$ of the compactification satisfies the following conditions:
$1\le k\le g+1$ and $k\equiv g+1~(\mod2)$ in the case $\ve=1$ and $0\le k\le g$ in the case $\ve=0$.
%These classification results were obtained by Weichold~\cite{Weichold}.
%In the case $\ve=1$ let $n=k$.
%In the case $\ve=0$ we choose $n\in\{k+1,\dots,g+1\}$ such that $n\equiv g+1~(\mod2)$.
%Let $\tilde g=(g+1-n)/2$.
Under these conditions the moduli space~$\calM_t$ has a natural topological structure, is connected and 
is the quotient $\calT_t/\Mod_t$,
where $\calT_t$ is homeomorphic to a real vector space of dimension $3g-3+3n_I+2m_I+2n_R+m_R$
and $\Mod_t$ is a discrete group.
\end{thm}

This result was previously only known for compact real algebraic curves~\cite{N1975,N1978moduli}.

\section{Sequential Sets of Automorphisms}

We will recall some standard results in complex analysis~\cite{Nbook2018}.
Consider the upper half plane $\hyp=\{z\in\c\st \Im(z)>0\}$.
Holomorphic and anti-holomorphic automorphisms of~$\hyp$
are isometries with respect to the hyperbolic metric $ds=\frac{|dz|}{\Im(z)}$ on~$\hyp$.
Geodesics in this geometry are half-circles with centre on the real axis and rays orthogonal to the real axis.
The group~$\Aut(\hyp)$ of bi-holomorphic automorphisms of~$\hyp$ is isomorphic to~$\PSL$
and consists of M\"obius transformations
\[C(z)=\frac{az+b}{cz+d},\quad\text{where}~a,b,c,d\in\r~\text{and}~ad-bc>0.\]
Bi-holomorphic automorphisms can be classified with respect to the fixed point behavior
of their action on~$\hyp$.
An element is called {\it hyperbolic} if it has two fixed points,
which lie on the boundary $\dd\hyp=\r\cup\{\infty\}$ of~$\hyp$.
A hyperbolic element with fixed points~$\al$, $\beta$ in~$\r$ is of the form
\[[\tau_{\al,\be}(\la)](z)=\frac{(\la\al-\be)z-(\la-1)\al\be}{(\la-1)z+(\al-\la\be)},\]
%\[
%  \tau_{\al,\be}(\la)
%  =\left[\frac{1}{(\al-\be)\cdot\sqrt{\la}}\cdot\begin{pmatrix}\la\al-\be&-(\la-1)\al\be\\ \la-1&\al-\la\be\end{pmatrix}\right],
%\]
where~$\la>0$.
A hyperbolic element with one fixed point at~$\infty$ is of the form
% comment: \lim\limits_{\al\to\infty}\tau_{\al,\be}
\[[\tau_{\infty,\be}(\la)](z)=\la z-(\la-1)\be\]
%\[\tau_{\infty,\be}(\la)=\left[\frac{1}{\sqrt{\la}}\cdot\begin{pmatrix}\la&-(\la-1)\be\\ 0&1\end{pmatrix}\right]\]
or
% comment: \lim\limits_{\be\to\infty}\tau_{\al,\be}
\[[\tau_{\al,\infty}(\la)](z)=\frac{1}{\la}z-\left(\frac{1}{\la}-1\right)\al,\]
%\[\tau_{\al,\infty}(\la)=\left[\frac{1}{\sqrt{\la}}\cdot\begin{pmatrix}1&(\la-1)\al\\ 0&\la\end{pmatrix}\right],\]
where $\al$ resp.\ $\be$ is the real fixed points and $\la>0$.
The parameter $\la>0$ is called the {\it shift parameter}.
The {\it axis} $\ell(g)$ of the element $g=\tau_{\al,\be}(\la)$ is the geodesic between the fixed points $\al$ and $\be$, oriented from $\be$ to $\al$ if $\la>1$ and from $\al$ to $\be$ if $\la<1$.
The element $g=\tau_{\al,\be}(\la)$ preserves the geodesic $\ell(g)$ and moves
the points on this geodesic in the direction of the orientation.
We call a hyperbolic element $\tau_{\al,\be}(\la)$ with $\la>1$ {\it positive} if $\al<\be$.
%The map $\la\mapsto\tau_{\al,\be}(\la)$ defines a homomorphism $\r_+\to G$ (with respect to the multiplicative structure on $\r_+$).
%We have
%\[(\tau_{\al,\be}(\la))^{-1}=\tau_{\al,\be}(\la^{-1})=\tau_{\be,\al}(\la).\]

\myskip
An element is called {\it parabolic} if it has one fixed point, which is on the boundary~$\dd\hyp$.
A parabolic element with a fixed point $\al$ is of the form
\[[\pi_{\al}(\la)](z)=\frac{(1-\la\al)z+\la\al^2}{-\la z+(1+\la\al)}.\]
%\[\pi_{\al}(\la)=\left[\begin{pmatrix}1-\la\al&\la\al^2\\ -\la&1+\la\al\end{pmatrix}\right].\]
A parabolic element with fixed point $\infty$ is of the form
\[[\pi_{\infty}(\la)](z)=z+\la.\]
%\[\pi_{\infty}(\la)=\left[\begin{pmatrix}1&\la\\ 0&1\end{pmatrix}\right].\]
We call a parabolic element $\pi_{\al}(\la)$ {\it positive} if $\la>0$.
%The map $\la\mapsto\pi_{\al}(\la)$ defines a homomorphism $\r\to G$ (with respect to the additive structure on $\r$).
%We have
%\[(\pi_{\al}(\la))^{-1}=\pi_{\al}(-\la).\]

\myskip
An element that is neither hyperbolic nor parabolic is called {\it elliptic}.
It has one fixed point that is in $\hyp$.
Given a base-point $x\in\hyp$ and a real number $\varphi$,
let $\rho_x(\varphi)$ denote the rotation through angle $\varphi$ about the point $x$.
Any elliptic element is of the form $\rho_x(\varphi)$, where $x$ is the fixed point.

\myskip
We will call hyperbolic and parabolic automorphisms of~$\hyp$ {\it shifts}.
Riemann surfaces are bi-holomorphic to quotients~$\hyp/\Ga$,
where $\Ga$ is a Fuchsian group that consists of shifts.
%without elliptic elements.
Sequential sets are special generating sets of such Fuchsian groups which were introduced in~\cite{N1972}.
They can be defined as follows:

\begin{mydef}
For two elements~$C_1$ and~$C_2$ in~$\Aut(\hyp)$ with finite fixed points in~$\r$ we say that $C_1<C_2$
if all fixed points of~$C_1$ are to the left of any fixed point of~$C_2$.
\end{mydef}

\begin{mydef}
\label{def-short-seq-set}
A triple of shifts $(C_1,C_2,C_3)$ in~$\Aut(\hyp)$ is a {\it sequential set}
%of type~$(0,n,m)$ with $n+m=3$
%the elements $C_i$ for~$i\le n$ are hyperbolic,
%the elements $C_i$ for~$i>n$ are parabolic,
if their product is
\[C_1\cdot C_2\cdot C_3=1,\]
and for some element $A\in\Aut(\hyp)$ the elements~$\{\tilde C_i=A C_i A^{-1}\}_{i=1,2,3}$ are positive, have finite fixed points and satisfy~$\tilde C_1<\tilde C_2<\tilde C_3$.
\end{mydef}

\begin{figure}[h]
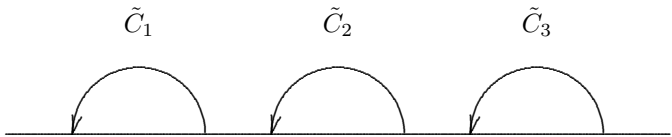

  \begin{center}
    \forcehmode
      \bgroup
        \beginpicture
          \setcoordinatesystem units <25 bp,25 bp>
          \multiput {\phantom{$\bullet$}} at -5 -1 5 2 /
          \circulararc 180 degrees from -2 0 center at -3 0
          \arrow <7pt> [0.2,0.5] from -3.99 0.04 to -4 0
          \put {$\tilde C_1$} [b] <0pt,\baselineskip> at -3 1
          \circulararc 180 degrees from  1 0 center at  0 0
          \arrow <7pt> [0.2,0.5] from -0.99 0.04 to -1 0
          \put {$\tilde C_2$} [b] <0pt,\baselineskip> at 0 1
          \circulararc 180 degrees from  4 0 center at  3 0
          \arrow <7pt> [0.2,0.5] from  2.01 0.04 to  2 0
          \put {$\tilde C_3$} [b] <0pt,\baselineskip> at 3 1
          \plot -5 0 5 0 /
        \endpicture
      \egroup
  \end{center}
  \caption{Sequential set}
%  \caption{Sequential set of type $(0,3,0)$}
  \label{fig-axes-seqset-0-3-0}
\end{figure}

Figure~\ref{fig-axes-seqset-0-3-0} illustrates the position of the axes of the elements~$\tilde C_i$
%for a sequential set of type $(0,3,0)$, \ie
when all elements are hyperbolic.
%It is clear what the similar picture looks like in presence of parabolic elements.

% FIGURE

\begin{mydef}
\label{def-long-seq-set}
A tuple of shifts $(C_1,\dots,C_r)$ in~$\Aut(\hyp)$ is a {\it sequential set}
%A {\it sequential set of type}~$(0,n,m)$ is an $(n+m)$-tuple of elements
%the elements $C_1,\dots,C_n$ are hyperbolic,
%the elements $C_{n+1},\dots,C_{n+m}$ are parabolic,
if for any $j\in\{2,\dots,r-1\}$ the triple $(C_1\cdots C_{j-1},C_j,C_{j+1}\cdots C_r)$
is a sequential set.
%(of type $(0,3,0)$, $(0,2,1)$, $(0,1,2)$ or~$(0,0,3)$).
\end{mydef}

\begin{mydef}
\label{def-genus-seq-set}
A {\it sequential set of type}~$(g,n,m)$ is a $(n+m+2g)$-tuple of shifts
\[(C_1,\dots,C_{n+m},A_1,B_1,\dots,A_g,B_g)\]
in~$\Aut(\hyp)$ such that
the elements $A_1,\dots,A_g,B_1,\dots,B_g$ and~$C_1,\dots,C_n$ are hyperbolic,
the elements $C_{n+1},\dots,C_{n+m}$ are parabolic,
and the tuple
\[(C_1,\dots,C_{n+m},A_1,B_1A_1^{-1}B_1^{-1},\dots,A_g,B_gA_g^{-1}B_g^{-1})\]
is a sequential set.
%of type $(0,2g+n,m)$.
\end{mydef}

\begin{figure}[h]
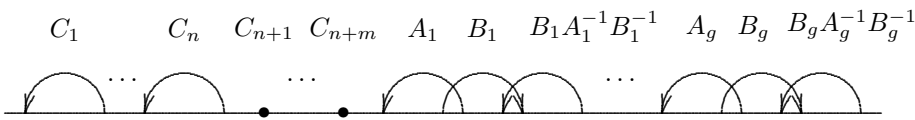

  \begin{center}
    \forcehmode
      \bgroup
        \beginpicture
          \setcoordinatesystem units <15 bp,15 bp>
          \multiput {\phantom{$\bullet$}} at -4.5 -1 17.5 2 /
          \circulararc 180 degrees from -2 0 center at -3 0
          \arrow <7pt> [0.2,0.5] from -3.99 0.04 to -4 0
          \put {$C_1$} [b] <0pt,\baselineskip> at -3 1
          \put {$\dots$} [b] <0pt,\baselineskip> at -1.5 0
          \circulararc 180 degrees from  1 0 center at  0 0
          \arrow <7pt> [0.2,0.5] from -0.99 0.04 to -1 0
          \put {$C_n$} [b] <0pt,\baselineskip> at 0 1
          \put {$\bullet$} at 2 0
          \put {$C_{n+1}$} [b] <0pt,\baselineskip> at 2 1
          \put {$\dots$} [b] <0pt,\baselineskip> at 3 0
          \put {$\bullet$} at 4 0
          \put {$C_{n+m}$} [b] <0pt,\baselineskip> at 4 1
          \circulararc 180 degrees from 7 0 center at 6 0
          \arrow <7pt> [0.2,0.5] from 5.01 0.04 to 5 0
          \put {$A_1$} [b] <0pt,\baselineskip> at 6 1
          \circulararc 180 degrees from 10 0 center at 9 0
          \arrow <7pt> [0.2,0.5] from 8.01 0.04 to 8 0
          \put {$B_1A_1^{-1}B_1^{-1}$} [bl] <-5pt,\baselineskip> at 9 1
          \circulararc 180 degrees from 8.5 0 center at 7.5 0
          \arrow <7pt> [0.2,0.5] from 8.49 0.04 to 8.5 0
          \put {$B_1$} [b] <0pt,\baselineskip> at 7.5 1
          \put {$\dots$} [b] <0pt,\baselineskip> at 11 0
          \circulararc 180 degrees from 14 0 center at 13 0
          \arrow <7pt> [0.2,0.5] from 12.01 0.04 to 12 0
          \put {$A_g$} [b] <0pt,\baselineskip> at 13 1
          \circulararc 180 degrees from 17 0 center at 16 0
          \arrow <7pt> [0.2,0.5] from 15.01 0.04 to 15 0
          \put {$B_gA_g^{-1}B_g^{-1}$} [bl] <-13pt,\baselineskip> at 16 1
          \circulararc 180 degrees from 15.5 0 center at 14.5 0
          \arrow <7pt> [0.2,0.5] from 15.49 0.04 to 15.5 0
          \put {$B_g$} [b] <-3pt,\baselineskip> at 14.5 1
          \plot -4.5 0 17.5 0 /
        \endpicture
      \egroup
  \end{center}
  \caption{Sequential set of type $(g,n,m)$}
  \label{fig-axes-seqset-g-n-m}
\end{figure}

Figure~\ref{fig-axes-seqset-g-n-m} illustrates the position of the axes and fixed points of the elements of a sequential set of type $(g,n,m)$.

% FIGURE

%\begin{figure}[H]
%\centering
%\includegraphics[width=1\textwidth]{FIG1-cropped}
%\caption{Sequential set of type $(g,n,m)$}
%\label{Figure1}
%\end{figure}

According to~\cite{N1978moduli,N1999b,Nbook}, a sequential set $\calV$ of type $(g,n,m)$
generates a Fuchsian group~$\Ga(\calV)$ that consists of shifts
%without elliptic elements
such that the surface $P=\hyp/\Ga$ is of genus~$g$ with $n$~holes and $m$~punctures.
We will say that $P$ is a surface of type~$(g,n,m)$.
%The isomorphism $\Phi:\lat\to\pi_1(P,p)$, induced by the natural projection $\Psi:\hyp\to P$,
%maps the sequential set $V$ to a standard basis of $\pi_1(P,p)$.

\section{Real Surfaces without Real Holes or Real Punctures}

\label{Section3}

Let the moduli space~$\calM_{g,n,m}$ be the space of classes of bi-holomorphic equivalence
of Riemann surfaces of type~$(g,n,m)$.
According to~\cite{N1978moduli,Nbook,Nbook2018}, any Riemann surface of type~$(g,n,m)$
is bi-holomorphic to a quotient~$\hyp/\Ga$, where $\Ga$ is a Fuchsian group
generated by a sequential set of type~$(g,n,m)$,
hence we can describe the moduli space~$\calM_{g,n,m}$
via the space of sequential sets of type~$(g,n,m)$.
Let $\tilde\calT_{g,n,m}$ be the set of all sequential sets of type~$(g,n,m)$.
The group~$\Aut(\hyp)\cong\PSL$ acts on the set~$\tilde\calT_{g,n,m}$ by conjugation.
Let $\calT_{g,n,m}=\tilde\calT_{g,n,m}/\PSL$ be the quotient of this action.
A sequential set of type~$(g,n,m)$ consists of $2g+n$~hyperbolic and $m$~parabolic automorphisms
which can be described by 
\[3(2g+n)+2m=6g+3n+2m\]
real parameters.
The relation
\[C_1\cdot\dots\cdot C_{n+m}\cdot[A_1,B_1]\cdot\dots\cdot[A_g,B_g]=1\]
implies some restrictions on these parameters, hence a point in $\tilde\calT_{g,n,m}$ 
can be described by $6g+3n+2m-3$ real parameters.
Taking into account the action of the $3$-dimensional group $\PSL$,
we can conclude that the space~$\calT_{g,n,m}$ is homeomorphic to $\r^{6g+3n+2m-6}$,
see~\cite{N1978moduli,Nbook,Nbook2018}.

\myskip
Sequential sets in the same orbit of~$\PSL$ correspond to bi-holomorphically equivalent Riemann surfaces,
hence the moduli space $\calM_{g,n,m}$ of Riemann surfaces of type~$(g,n,m)$
up to bi-holomorphic equivalence is of the form
\[\calM_{g,n,m}=\calT_{g,n,m}/\Mod_{g,n,m},\]
where $\Mod_{g,n,m}$ is a discrete group of automorphisms of a surface of type~$(g,n,m)$
that acts on~$\calT_{g,n,m}$.
Therefore
\[\calM_{g,n,m}\simeq\r^{6g-6+3n+2m}/\Mod_{g,n,m}.\]
This fact is known as Fricke-Klein Theorem~\cite{FK}.

\begin{figure}[h]
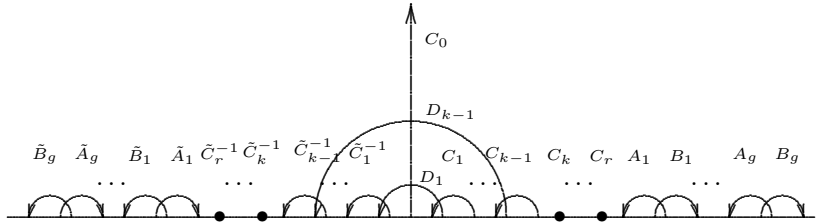

  \begin{center}
    \forcehmode
      \bgroup
        \beginpicture
          \setcoordinatesystem units <8 bp,8 bp>
          \multiput {\phantom{$\bullet$}} at -18 -1 20 2 1 11 /
          \circulararc 180 degrees from -15 0 center at -16 0
          \arrow <5pt> [0.2,0.5] from -16.99 0.04 to -17 0
          \put {\tiny{$\tilde B_g$}} [b] <-2pt,\baselineskip> at -16 1
          \circulararc 180 degrees from -13.5 0 center at -14.5 0
          \arrow <5pt> [0.2,0.5] from -13.51 0.04 to -13.5 0
          \put {\tiny{$\tilde A_g$}} [b] <2pt,\baselineskip> at -14.5 1
          \put {$\dots$} [b] <0pt,\baselineskip> at -13 0
          \circulararc 180 degrees from -10.5 0 center at -11.5 0
          \arrow <5pt> [0.2,0.5] from -12.49 0.04 to -12.5 0
          \put {\tiny{$\tilde B_1$}} [b] <-2pt,\baselineskip> at -11.5 1
          \circulararc 180 degrees from -9 0 center at -10 0
          \arrow <5pt> [0.2,0.5] from -9.01 0.04 to -9 0
          \put {\tiny{$\tilde A_1$}} [b] <2pt,\baselineskip> at -10 1
          \put {$\bullet$} at -8 0
          \put {\tiny{$\tilde C_r^{-1}$}} [b] <0pt,\baselineskip> at -8 1
          \put {$\dots$} [b] <0pt,\baselineskip> at -7 0
          \put {$\bullet$} at -6 0
          \put {\tiny{$\tilde C_k^{-1}$}} [b] <0pt,\baselineskip> at -6 1
          \circulararc 180 degrees from -3 0 center at -4 0
          \arrow <5pt> [0.2,0.5] from -4.99 0.04 to -5 0
          \put {\tiny{$\tilde C_{k-1}^{-1}$}} [b] <5pt,\baselineskip> at -4 1
          \put {$\dots$} [b] <0pt,\baselineskip> at -2.5 0
          \circulararc 180 degrees from  0 0 center at  -1 0
          \arrow <5pt> [0.2,0.5] from -1.99 0.04 to -2 0
          \put {\tiny{$\tilde C_1^{-1}$}} [b] <0pt,\baselineskip> at -1 1
          \arrow <7pt> [0.2,0.5] from 1 0 to 1 10
          \put {\tiny{$C_0$}} [bl] <5pt,0pt> at 1 8
          \circulararc 180 degrees from 4 0 center at 3 0
          \arrow <5pt> [0.2,0.5] from 2.01 0.04 to 2 0
          \put {\tiny{$C_1$}} [b] <0pt,\baselineskip> at 3 1
          \put {$\dots$} [b] <0pt,\baselineskip> at 4.5 0
          \circulararc 180 degrees from 7 0 center at  6 0
          \arrow <5pt> [0.2,0.5] from 5.01 0.04 to 5 0
          \put {\tiny{$C_{k-1}$}} [b] <-3pt,\baselineskip> at 6 1
          \put {$\bullet$} at 8 0
          \put {\tiny{$C_k$}} [b] <0pt,\baselineskip> at 8 1
          \put {$\dots$} [b] <0pt,\baselineskip> at 9 0
          \put {$\bullet$} at 10 0
          \put {\tiny{$C_r$}} [b] <0pt,\baselineskip> at 10 1
          \circulararc 180 degrees from 13 0 center at 12 0
          \arrow <5pt> [0.2,0.5] from 11.01 0.04 to 11 0
          \put {\tiny{$A_1$}} [b] <-2pt,\baselineskip> at 12 1
          \circulararc 180 degrees from 14.5 0 center at 13.5 0
          \arrow <5pt> [0.2,0.5] from 14.49 0.04 to 14.5 0
          \put {\tiny{$B_1$}} [b] <2pt,\baselineskip> at 13.5 1
          \put {$\dots$} [b] <0pt,\baselineskip> at 15 0
          \circulararc 180 degrees from 18 0 center at 17 0
          \arrow <5pt> [0.2,0.5] from 16.01 0.04 to 16 0
          \put {\tiny{$A_g$}} [b] <-2pt,\baselineskip> at 17 1
          \circulararc 180 degrees from 19.5 0 center at 18.5 0
          \arrow <5pt> [0.2,0.5] from 19.49 0.04 to 19.5 0
          \put {\tiny{$B_g$}} [b] <2pt,\baselineskip> at 18.5 1
          \circulararc 180 degrees from 2.5 0 center at 1 0
          \arrow <5pt> [0.2,0.5] from -0.49 0.04 to -0.5 0
          \put {\tiny{$D_1$}} [b] <8pt,0pt> at 1 1.5
          \circulararc 180 degrees from 5.5 0 center at 1 0
          \arrow <5pt> [0.2,0.5] from -3.49 0.04 to -3.5 0
          \put {\tiny{$D_{k-1}$}} [b] <15pt,0pt> at 1 4.5
          \plot -18 0 20 0 /
        \endpicture
      \egroup
  \end{center}
  \caption{Sequential set for a real curve of type $(g,k,\ve|2n_I,2m_I)$}
  \label{fig-axes-seqset-type-g-k-ve-2nI-2mI}
\end{figure}

\myskip
Now let us consider real algebraic curves.
Recall that the topological type of a real curve without real holes and real punctures
is determined by the type~$(g,k,\ve)$ of the corresponding compactified real curve
and the numbers~$2n_I$ and~$2m_I$ of non-real holes and punctures.
A slight modification of the proofs in~\cite{N1978moduli,N1999a,N1999b} leads to the following description
of the corresponding Fuchsian groups.
A real curve of type~$t=(g,k,\ve|2n _I,2m_I)$ can be constructed using a sequential set
\[\calV=(C_0,C_1,\dots,C_r,A_1,\dots,A_h,B_1,\dots,B_h)\]
of type~$(h,g-2h+1+n_I,m_I)$ and hyperbolic automorphisms
\[\calD=(D_1,\dots,D_{g-2h}),\]
where the number~$h$ and the hyperbolic automorphisms~$\calD$ depend on the type~$t$.
The hyperbolic automorphisms~$\calD$ are constructed using reflections in geodesics in the hyperbolic space~$\hyp$.
Let $R_j$ denote the reflection in the axis of the hyperbolic automorphism~$C_j$.
%and let $\tilde c=\bar c\sqrt{c}$.
For~$\ve=1$, let $h=(g-k+1)/2$ and $D_j=R_0R_j$ for~$j=1,\dots,g-2h=k-1$.
For a hyperbolic automorphism~$C$,
let $\sqrt{C}$ be the hyperbolic automorphism such that~$(\sqrt{C})^2=C$.
For~$\ve=0$ and~$k>0$, let $h=0$, $D_j=R_0R_j$ for $j=1,\dots,k-1$
and $D_j=R_0R_j\sqrt{C_j}$ for~$j=k,\dots,g$.
For~$\ve=0$ and~$k=0$, let $h=0$ and $D_j=\sqrt{C_0}R_0R_j\sqrt{C_j}$ for $j=1,\dots,g$.
Figure~\ref{fig-axes-seqset-type-g-k-ve-2nI-2mI} shows the axes of the elements
with the notation $\tilde A_i=R_0 A_i R_0$, $\tilde B_i=R_0 B_i R_0$, $C_j=R_0 C_j R_0$.

% FIGURE

%\begin{figure}[H]
%\centering
%\includegraphics[width=1\textwidth]{FIG2-cropped}
%\caption{XXX.}
%\label{Figure2}
%\end{figure}

We can show that automorphisms in $\calV\cup\calD$ generate a Fuchsian group~$\Ga=\Ga(\calV\cup\calD)$
such that the Riemann surface $P=\hyp/\Ga$ is of type~$(g,2n_I,2m_I)$.
Let $\si=R_0$ for~$k>0$ and $\si=\sqrt{C_0} R_0$ for~$k=0$, then $\si \Ga \si=\Ga$.
It follows that $\si$ generates a real curve~$(P,\tau)$ of type~$(g,k,\ve|2n_I,2m_I)$.
Moreover, it can be shown that every real curve of type~$t$ can be constructed in this way.

\myskip
The hyperbolic metric on the upper half-plane~$\hyp$
induces a hyperbolic metric on the quotient~$P=\hyp/\Ga$.
This metric reflects the geometry of the Fuchsian group~$\Ga$.
For instance, among all simple closed curves around the hole that corresponds to the generator~$C_j$
there is a unique shortest curve~$c_j$.
This curve is the image of the axis of~$C_j$ under the natural projection~$\hyp\to P$.

\myskip
Recall that real curves of type~$t=(g,k,\ve|2n_I,2m_I)$
are generated by sequential sets of type~$(h,g-2h+n_I,m_I)$,
whereby conjugate sequential sets correspond to bi-holomorphically equivalent real curves.
The space~$\calT_{h,g-2h+1+n_I,m_I}$ of conjugacy classes of such sequential sets is homeomorphic to
a real vector space of dimension
\[6h+3(g-2h+1+n_I)+2m_I-6=3g-3+3n_I+2m_I.\]
The moduli space~$\calM_{(g,k,\ve | 2n_I,2m_I)}$ is obtained as the quotient of $\calT_{h,g-2h+1+n_I,m_I}$
by the discrete group~$\Mod_t$ of homotopy classes of those automorphisms that commute with~$\tau$.
Thus $\calM_{(g,k,\ve | 2n_I,2m_I)}$ is homeomorphic to
\[\r^{3g-3+3n_I+2m_I}/\Mod_t.\]

\section{Real Surfaces of Genus Zero}

\label{Section4}

\begin{figure}[h]
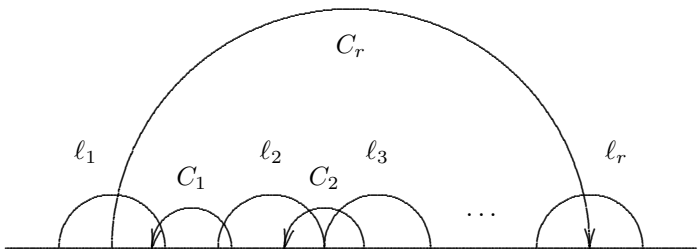

  \begin{center}
    \forcehmode
      \bgroup
        \beginpicture
          \setcoordinatesystem units <20 bp,20 bp>
          \multiput {\phantom{$\bullet$}} at -1 -1 12 5 /
          \circulararc 180 degrees from 2 0 center at 1 0
          \put {$\ell_1$} [b] <0pt,\baselineskip> at 0.5 1
          \circulararc 180 degrees from 5 0 center at 4 0
          \put {$\ell_2$} [b] <0pt,\baselineskip> at 4 1
          \circulararc 180 degrees from 7 0 center at 6 0
          \put {$\ell_3$} [b] <0pt,\baselineskip> at 6 1
          \put {$\dots$} [b] <0pt,\baselineskip> at 8 0
          \circulararc 180 degrees from 11 0 center at 10 0
          \put {$\ell_r$} [b] <0pt,\baselineskip> at 10.5 1
          \circulararc 180 degrees from 3.25 0 center at 2.5 0
          \arrow <7pt> [0.2,0.5] from 1.76 0.04 to 1.75 0
          \put {$C_1$} [b] <0pt,\baselineskip> at 2.5 0.5
          \circulararc 180 degrees from 5.75 0 center at 5 0
          \arrow <7pt> [0.2,0.5] from 4.26 0.04 to 4.25 0
          \put {$C_2$} [b] <0pt,\baselineskip> at 5 0.5
          \circulararc 180 degrees from 10 0 center at 5.5 0
          \arrow <7pt> [0.2,0.5] from 10 0.04 to 10 0
          \put {$C_r$} [b] <0pt,\baselineskip> at 5.5 3
          \plot -1 0 12 0 /
        \endpicture
      \egroup
  \end{center}
  \caption{Real curve of genus zero}
  \label{fig-genus-zero}
\end{figure}

We will now prove Theorem~\ref{thm-moduli} for real surfaces of genus zero without non-real holes or punctures.
In this case the compactified real curve has exactly one oval and the number of holes and punctures on this oval determines the type~$t$ of the real surface.
Let us consider a set $\ell_1,\dots,\ell_r$, $r>2$, of pairwise disjoint geodesics in the hyperbolic plane~$\hyp$ such that the end points of~$\ell_i$ are to the left of the endpoints of~$\ell_{i+1}$.
Let $R_i$ be the reflection in~$\ell_i$.
Product of two hyperbolic reflections
is parabolic if their axes share exactly one endpoint
and hyperbolic if the closures of their axes are disjoint.
We can choose the geodesics~$\ell_i$ in such a way that the distribution of parabolic and hyperbolic elements
among the products $C_r=R_rR_1,C_1=R_1R_2,C_2=R_2R_3,\dots,C_{r-1}=R_{r-1}R_r$ corresponds to the distribution of punctures and holes in the topological type~$t$, see Figure~\ref{fig-genus-zero}.

% FIGURE

%\begin{figure}[H]
%\centering
%\includegraphics[width=1\textwidth]{FIG3-cropped}
%\caption{XXX.}
%\label{Figure3}
%\end{figure}

\noindent
Consider the shifts
\[F_1=R_1R_i=C_1\cdots C_{i-1}\quad\text{and}\quad F_2=R_{i+1}R_r=C_{i+1}\cdots C_{r-1}.\]
Then $(F_1,C_i,F_2)$ is a sequential set and hence so is $(C_1,\dots,C_r)$.
Therefore the shifts~$C_1,\dots,C_r$ generate a Fuchsian group~$\Ga$
%that consists of shifts
%without elliptic elements
and $P_r=\hyp/\Ga$ is a Riemann surface of genus zero.
Moreover, $R_i\Ga R_i=\Ga$,
hence the reflections~$R_i$ induce an anti-holomorphic involution $\tau:P_r\to P_r$
such that all holes and punctures are real with respect to~$\tau$
and their distribution corresponds to the distribution prescribed in the topological type~$t$.

\myskip
We will now prove that this construction gives all real curves~$(P_r,\tau)$ of genus zero
with only real holes and punctures.
To this end, given such a real curve, we want to identify the corresponding geodesics~$\ell_1,\dots,\ell_r$ in~$\hyp$.
Consider a connected component of~$P_r\backslash P_r^{\tau}$.
Let $Q_r$ be its pre-image under the uniformisation map $\hyp\to P_r$.
The boundary of~$Q_r$ consists of geodesic segments in~$\hyp$.
We can extend these geodesic segments to geodesics~$\ell_i$.

\myskip
Sets of geodesics are determined by their endpoints.
The number of endpoints is $2n+m$, where $n$ and $m$ are the numbers of holes and punctures in the type~$t$ respectively.
Taking into account the action of the $3$-dimensional group $\PSL$ of automorphisms of~$\hyp$ 
on the set of sets of geodesics,
we can conclude that the set of orbits of this action is homeomorphic to~$\r^{2n+m-3}$,
thus $\calM_t\simeq\r^{2n+m-3}/\Mod$,
where $\Mod$ is a discrete group determined by the numbering of the holes.

\begin{figure}[h]
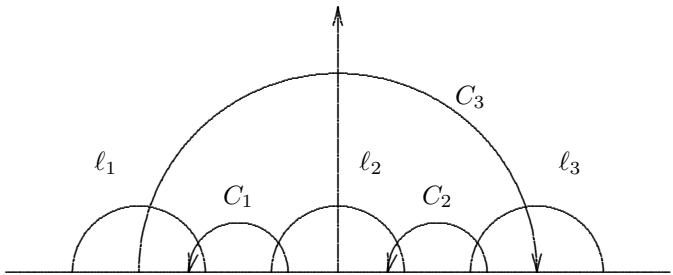

  \begin{center}
    \forcehmode
      \bgroup
        \beginpicture
          \setcoordinatesystem units <25 bp,25 bp>
          \multiput {\phantom{$\bullet$}} at -1 -1 9 4 /
          \circulararc 180 degrees from 2 0 center at 1 0
          \put {$\ell_1$} [b] <0pt,\baselineskip> at 0,5 1
          \circulararc 180 degrees from 5 0 center at 4 0
          \put {$\ell_2$} [b] <0pt,\baselineskip> at 4.5 1
          \circulararc 180 degrees from 8 0 center at 7 0
          \put {$\ell_3$} [b] <0pt,\baselineskip> at 7.5 1
          \circulararc 180 degrees from 3.25 0 center at 2.5 0
          \arrow <7pt> [0.2,0.5] from 1.76 0.04 to 1.75 0
          \put {$C_1$} [b] <0pt,\baselineskip> at 2.5 0.5
          \circulararc 180 degrees from 6.25 0 center at 5.5 0
          \arrow <7pt> [0.2,0.5] from 4.76 0.04 to 4.75 0
          \put {$C_2$} [b] <0pt,\baselineskip> at 5.5 0.5
          \circulararc 180 degrees from 7 0 center at 4 0
          \arrow <7pt> [0.2,0.5] from 7 0.04 to 7 0
          \put {$C_3$} [b] <0pt,\baselineskip> at 6 2
          \arrow <7pt> [0.2,0.5] from 4 0 to 4 4
          \plot -1 0 9 0 /
        \endpicture
      \egroup
  \end{center}
  \caption{Real curve of genus zero with $3$~holes}
  \label{fig-hexagon}
\end{figure}

\myskip
We will now construct the space~$\calM^*$ of real curves of genus zero
with $3$~holes of which exactly one is real.
Let us revisit the previous construction with $r=3$.
We consider a set $\ell_1,\ell_2,\ell_3$ of pairwise disjoint geodesics in the hyperbolic plane~$\hyp$
such that the end points of~$\ell_i$ are to the left of the endpoints of~$\ell_{i+1}$.
Let $R_i$ be the reflection in~$\ell_i$.
We can choose the geodesics~$\ell_i$ in such a way that all products
$C_3=R_3R_1,C_1=R_1R_2,C_2=R_2R_3$ are hyperbolic.
Let $c_i$ be the axis of $C_i$ for $i=1,2,3$.
In this case the geodesics $\ell_1,\ell_2,\ell_3$ and $c_1,c_2,c_3$ bound a right-angled hexagon~$Q_3$ in~$\hyp$.
The type of the hexagon~$Q_3$ is determined up to bi-holomorphic equivalence
by the lengths of its sides that are contained in the geodesics $c_1,c_2,c_3$, see~\cite{Thurston}.
These lengths are given by $(\la_1/2,\la_2/2,\la_3/2)$,
where $\la_i=\la(C_i)$ is the shift parameter of the hyperbolic element~$C_i$.
They are completely determined by the lengths of the minimal geodesics
around the holes of the surface~$P_3$.
If~$\la_1=\la_2$ then the hexagon~$Q_3$ is symmetric
with respect to the geodesic orthogonal to~$c_3$, see Figure~\ref{fig-hexagon}.
The reflection in this geodesic generates an anti-holomorphic involution $\si:P_3\to P_3$
that interchanges the holes~$c_1$ and~$c_2$ and maps the hole~$c_3$ to itself.
Note that $\la_1$ and~$\la_3$ completely determine the real curve~$(P_3,\si)$
and hence the space~$\calM^*$ of such real curves is homeomorphic to~$\r^2$.

% FIGURE

%\begin{figure}[H]
%\centering
%\includegraphics[width=1\textwidth]{FIG4-cropped}
%\caption{Real curve of genus zero with $3$~holes}
%\label{fig-hexagon}
%\end{figure}

\myskip
Now let us construct the space~$\calM^*_s$ of real curves $(P,\tau)$ of genus zero with two non-real holes
and $r$~real holes and punctures, generating an oval of type~$s$.
Consider the minimal closed geodesic~$c_0$ that separates non-real holes from real holes and punctures.
The geodesic~$c_0$ divides the real curve~$(P,\tau)$ into two real curves $(P_{r+1},\tau)$ and~$(P_3,\si)$
of the types considered above.
On the other hand, starting with the real curves $(P_{r+1},\tau)$ and~$(P_3,\si)$, we can glue them together to form a real curve of type~$t$
if the length of the minimal closed geodesic around the real hole on~$(P_3,\si)$
is equal to the length of the minimal closed geodesic around the hole~$C_1$ on~$(P_{r+1},\tau)$.
Therefore $M_t\simeq\r^{2n_R+m_R-3}/\Mod$, where $n_R$ and~$m_R$ are the numbers of real holes and real punctures in the type~$s$ respectively.

\section{Real Surfaces of Any Type}

\label{Section5}

We will now consider the space of real surfaces of any type.
Such real surfaces can be obtained by identifying pairs of non-real symmetric holes on the kinds of surfaces considered in sections~\ref{Section3} and~\ref{Section4}.

\myskip
Consider a real curve~$(P,\tau)$ of type~$t$ without real holes or real punctures.
Let $c$ be a minimal closed geodesic around a non-real hole on~$(P,\tau)$.
Let~$\rho$ be the length of~$c$.
%Let $x_i$ be the point constructed in section~\ref{Section3}.
Let $c_0$ be as described in section~\ref{Section3}.
Among all geodesic segments connecting the curves~$c$ and~$c_0$ there is a unique shortest segment.
Let $x$ be the end point of this segment on the curve~$c$.

\myskip
Consider a real curve~$(Q,\si)$ of type~$s$ with exactly two non-real holes and without non-real punctures.
Let $u$ and~$\si(u)$ be the minimal simple closed geodesics around the non-real holes on~$(Q,\si)$.
We can choose $(Q,\si)$ in such a way that the length of~$u$ and~$\si(u)$ is equal to~$\rho$.
%Let $y_i$ be the point constructed in section~\ref{Section4}.
Let $c_0$ be the minimal closed geodesic on $(Q,\si)$
that separates non-real holes from real holes and real punctures.
Among all geodesic segments connecting the curves~$u$ and~$c_0$
there is a unique shortest segment.
Let $y$ be the end point of this segment on the curve~$u$.   

\myskip
The family of isometries~$\varphi$ mapping the geodesic~$c$ to the geodesic~$u$ is real one-dimensional and can be parametrised by the distance between the points~$\varphi(x)$ and~$y$.
We can identify the holes~$\tau(c)$ and~$\si(u)$ via the isometry~$\si\circ\varphi\circ\tau$.
Thus we can glue a real curve~$(Q,\si)$ of genus zero and type~$s$
into a non-real hole of a real surface~$(P,\tau)$ of type~$t$.
Repeatedly gluing in real surfaces of genus zero we can obtain a real curve of any type.

\myskip
Consider the space~$\calM_{(t|s_1,\dots,s_r)}$ 
of real curves of topological type~$(t|s_1,\dots,s_r)$,
where $t$ is the topological type of a real curve without holes and punctures
and $s_i$ are the types of ovals.
The results of sections~\ref{Section4} and~\ref{Section5} imply that
\[\calM_{(t | s_1,\dots,s_r)}=\calM_t\times\r^{2n_R+m_R}.\]
Recall that we have shown in section~\ref{Section3} that
\[\calM_t\simeq\r^{3g-3+3n_I+2m_I}/\Mod_t,\]
where~$t=(g,k,\ve|2n_I,2m_I)$.

%%\nocite{*}
 %%\bibliographystyle{amsalpha}
%%\bibliography{realforms}

\def\cprime{$'$}
\providecommand{\bysame}{\leavevmode\hbox to3em{\hrulefill}\thinspace}
\providecommand{\MR}{\relax\ifhmode\unskip\space\fi MR }
\providecommand{\MRhref}[2]{%
  \href{http://www.ams.org/mathscinet-getitem?mr=#1}{#2}
}
\providecommand{\href}[2]{#2}

\end{document}